\documentclass[11pt,reqno]{amsart}
\usepackage{amsmath, amsthm, amssymb}
\usepackage{lipsum,mathtools,graphicx}
\usepackage{xspace,tikz-cd}
\usepackage[smalltableaux]{ytableau}
\newtheorem{thm}{Theorem}
\newtheorem{lemma}{Lemma}
\newtheorem{corollary}{Corollary}
\newtheorem{proposition}{Proposition}

\theoremstyle{definition}
\newtheorem*{remarks}{Remarks}

\newcommand{\qw}{$q$-Whittaker\xspace}
\newcommand{\mhl}{modified Hall-Littlewood\xspace}
\newcommand{\wl}[1][q]{W_\lambda(X_n; #1)}
\newcommand{\qpl}[1][\lambda]{Q'_{#1}(X_n; q)}
\newcommand{\htl}{\widetilde{H}_{\lambda}(X_n; q,t)}
\newcommand{\naturals}{\mathbb{N}}
\newcommand{\scf}{\mathcal{F}(\lambda)}
\newcommand{\csf}[1][\lambda]{\mathrm{CSF}(#1)}
\newcommand{\wdf}{\mathrm{WDF}(\lambda)}
\newcommand{\GT}{\mathrm{GT}(\lambda)}
\newcommand{\qbinom}[2]{{#1 \brack #2}_q}
\newcommand{\pop}[1][\lambda]{\mathrm{POP}(#1)}
\newcommand{\popp}[1][\lambda]{\mathrm{POPP}(#1)}
\newcommand{\barsig}{\overline{\sigma}}
\newcommand{\bartau}{\overline{\tau}}
\newcommand{\fdag}{F^\dag}
\newcommand{\bgam}{\overline{\gamma}}
\newcommand{\rootthree}{1.7320508}

\DeclareMathOperator{\len}{len}
\DeclareMathOperator{\inv}{inv}
\DeclareMathOperator{\refinv}{refinv}
\DeclareMathOperator{\quinv}{quinv}
\DeclareMathOperator{\maj}{maj}
\DeclareMathOperator{\dg}{dg}
\DeclareMathOperator{\dgh}{\widehat{\dg}}
\DeclareMathOperator{\NE}{NE}
\DeclareMathOperator{\SE}{SE}
\DeclareMathOperator{\br}{br}
\DeclareMathOperator{\pr}{pr}
\DeclareMathOperator{\bcomp}{boxcomp}
\DeclareMathOperator{\area}{area}
\DeclareMathOperator{\rsort}{rsort}

\DeclareMathOperator{\SSYT}{SSYT}
\DeclareMathOperator{\wt}{wt}
\DeclareMathOperator{\zcount}{zcount}
\DeclareMathOperator{\zcb}{\overline{\zcount}}
\DeclareMathOperator{\coarm}{coarm}
\DeclareMathOperator{\arm}{arm}
\DeclareMathOperator{\Inv}{Inv}
\DeclareMathOperator{\Des}{Des}
\DeclareMathOperator{\cells}{cells}
\DeclareMathOperator{\splice}{splice}
\DeclareMathOperator{\dsplice}{dsplice}
\definecolor{cadmiumgreen}{rgb}{0.0, 0.42, 0.24}
\definecolor{chocolate}{rgb}{0.48, 0.25, 0.0}
\definecolor{darkpastelgreen}{rgb}{0.01, 0.75, 0.24}
\definecolor{cadmiumred}{rgb}{0.89, 0.0, 0.13}
\definecolor{orangered}{rgb}{1.0, 0.27, 0.0}
\definecolor{light-gray}{gray}{0.9} 

\usepackage[backend=bibtex]{biblatex}
\begin{document}
\title[Monomial expansions]{Monomial expansions  for
$q$-Whittaker and modified Hall-Littlewood polynomials (extended abstract)}

\author[Bhattacharya]{Aritra Bhattacharya}
\address{The Institute of Mathematical Sciences, A CI of Homi Bhabha National Institute, Chennai 600113, India}
\email{baritra@imsc.res.in}.

\author[Ratheesh]{T V Ratheesh}
\address{The Institute of Mathematical Sciences, A CI of Homi Bhabha National Institute, Chennai 600113, India}
\email{ratheeshtv@imsc.res.in}

\author[Viswanath]{Sankaran Viswanath}
\address{The Institute of Mathematical Sciences, A CI of Homi Bhabha National Institute, Chennai 600113, India}
\email{svis@imsc.res.in}

\thanks{The authors acknowledge partial support under a DAE Apex Grant to the Institute of Mathematical Sciences.}

\begin{abstract}
We consider the monomial expansion of the $q$-Whittaker polynomials given by the fermionic formula and via the {\em inv} and {\em quinv} statistics. We construct bijections between the parametrizing sets of these three models which preserve the $x$- and $q$-weights, and which are compatible with natural projection and branching maps. We apply this to the limit construction of local Weyl modules and obtain a new character formula for the basic representation of $\widehat{\mathfrak{sl}_n}$. Finally, we indicate how our main results generalize to the modified Hall-Littlewood case.
\end{abstract}
\maketitle

\section{Introduction}
Let $\lambda$ be a partition. For $n \geq 1$, let $X_n$ denote the tuple of indeterminates $x_1, x_2, \cdots, x_n$.
The \qw polynomial $\wl$ and the \mhl polynomial $\qpl$ are well-studied specializations of the modified Macdonald polynomial. Several different monomial expansions for these polynomials are known. In this article, our focus will be on three of these: the so-called {\em fermionic formulas} \cite[(0.2), (0.3)]{Kirillov_newformula} and the inv- and quinv-expansions arising from specializations of the formulas of Haglund-Haiman-Loehr \cite{HHL-I} and  Ayyer-Mandelshtam-Martin \cite{AMM}.

We recall that the Schur expansion of the $\wl$ (resp. $\qpl$) has certain $q$-Kostka polynomials as coefficients \cite[]{Kirillov_newformula}. In turn, this implies yet another monomial expansion, with the underlying indexing set involving pairs of semistandard Young tableaux of conjugate (resp. equal) shapes. This relates to the inv-expansion via the RSK correspondence \cite{HHL-I}.

The fermionic formula, expressed as a sum of products of $q$-binomials, is seemingly of a very different nature from all the other monomial expansions, and should probably viewed as a kind of compression of these formulas. Recently, Garbali-Wheeler \cite{GarbaliWheeler} obtained a general formula of the fermionic kind for the full modified Macdonald polynomial $\htl$.

The purpose of this article is to bijectively reconcile the fermionic formula with both the inv- and quinv-expansions. We construct bijections between the underlying sets of these three models which (i) preserve the $x$- and $q$-weights, and (ii) are compatible with natural projection and branching maps.

As a corollary, we obtain bijections between the inv- and quinv-models in the \qw and \mhl specializations, partially answering a question of \cite{AMM}. We find that the $inv$- and $quinv$-models are related by the simple {\em box-complementation} map of the fermionic model and that $inv+quinv$ is a constant on fibers of the natural projection. We also apply this to the limit construction for Weyl modules \cite{fourier-littelmann,RRV-CLpop} and 
obtain an apparently new character formula for the basic representation of the affine Lie algeba $\widehat{\mathfrak{sl}_n}$.

In this extended abstract, we describe the \qw polynomials in greater detail, contenting ourselves with brief remarks about the \mhl case in \S\ref{sec:concl-rem} due to space limitations. Complete proofs will appear in \cite{brv-fullversion}.

\section{Specializations of $\htl$}
Given a partition $\lambda=(\lambda_1 \geq \lambda_2 \geq \cdots)$, we will draw its Young diagram $\dg(\lambda)$ following the English convention, as a left-up justified array of boxes, with $\lambda_i$ boxes in the $i$th row from the top. The boxes are called the cells of $\dg(\lambda)$. We let $|\lambda|:=\sum_i \lambda_i$.
Fix $n \geq 1$ and let $\scf$ denote the set of all maps (``fillings'') $F: \dg(\lambda) \to [n]$ where $[n]=\{1, 2, \cdots, n\}$. If the values of $F$ strictly increase (resp. weakly decrease) as we move down a column, we say $F$ is a {\em column strict filling} (CSF) (resp. {\em weakly decreasing filling} (WDF)\footnote{These latter ones may be easily transformed into the familiar {\em tabloids} by transposing rows and columns and replacing $i \mapsto n-i+1$}), and denote the set of such fillings by $\csf$ (resp. $\wdf$).
The $x$-weight of a filling $F$ is the monomial $x^F:=\displaystyle\prod_{c \in \dg(\lambda)} x_{F(c)}$.

We recall that the modified Macdonald polynomial $\htl$ is a symmetric polynomial in the $x_i$ with $\naturals[q,t]$ coefficients. We expand this in powers of $t$; our interest lies in the coefficients of the lowest and highest powers \cite[(3.1)]{bergeron2020survey}:
\begin{equation}\label{eq:berg1}
  \htl = \mathcal{H}_\lambda(X_n; q) t^0 + \cdots + \wl t^{\eta(\lambda)}
  \end{equation}
where $\eta(\lambda) = \sum_{j \geq 1} \binom{\lambda'_j}{2}$ where $\lambda'_j$ denote the parts of the partition conjugate  to $\lambda$. The $\wl$ is the \qw polynomial. The $q$-reversal (or reciprocal) polynomial of $\mathcal{H}_\lambda(X_n; q)$ coincides with the \mhl polynomial $\qpl[\lambda']$ where $\lambda'$ is the partition conjugate to $\lambda$, i.e., $q^{\eta(\lambda')} \mathcal{H}_\lambda(X_n; q^{-1}) = \qpl[\lambda']$. These are further related to each other by $\omega\wl = \qpl[\lambda']$ where $\omega$ is the classical involution on the ring of symmetric polynomials.

Following Haglund-Haiman-Loehr \cite{HHL-I} and Ayyer-Mandelshtam-Martin \cite{AMM}, there are statistics {\em inv, quinv} and {\em maj} on $\scf$ such that
\begin{equation}\label{eq:invquinvmaj}
\htl = \sum_{F \in \scf} x^F q^{v(F)} t^{\maj(F)}
\end{equation}
where $v \in \{\inv, \quinv\}$.
The next lemma follows directly from the definition of $\maj$ \cite{HHL-I}:
\begin{lemma}\label{lem:majminmax}
  Let $F \in \scf$. Then (i) $\maj(F) = \eta(\lambda)$ iff $F \in \csf$, and (ii) $\maj(F)=0$ iff $F \in \wdf$.
\end{lemma}

\noindent
Putting together \eqref{eq:berg1}, \eqref{eq:invquinvmaj} and Lemma~\ref{lem:majminmax}, we obtain for $v \in \{\inv,\quinv\}$:
\begin{align}
  \wl &= \sum_{F \in \csf} x^F q^{v(F)} \label{eq:wlinvquinv}\\
  \qpl[\lambda'] &= \sum_{F \in \wdf} x^F q^{\eta(\lambda') - v(F)} \label{eq:qlpinvquinv}
\end{align}
These are in fact symmetric in the $x$-variables  and can be viewed as  expansions in terms of the monomial symmetric functions in $x_1, x_2, \cdots, x_n$.

\section{Fermionic formula for $\wl$}\label{sec:fermfor}

Let $n \geq 1$ and $\lambda = (\lambda_1 \geq \lambda_2 \geq \cdots \geq \lambda_n \geq 0)$ be a partition with at most $n$ nonzero parts. Let $\GT$ denote the set of integral Gelfand-Tsetlin (GT) patterns with bounding row $\lambda$. Given $T \in \GT$, we denote its entries by $T^j_i$ for $1 \leq i \leq j \leq n$ as in Figure~\ref{fig:gt-ssyt}. It will also be convenient to define $T^j_{j+1}=0$ for all $1 \leq j \leq n$.
We define the North-East and South-East differences of $T$ by:
$\NE_{ij}(T) = T^{j+1}_i - T^j_i$ and $\SE_{ij}(T)= T^j_i - T^{j+1}_{i+1}$
for $1 \leq i \leq (j+1) \leq n$. The GT inequalities ensure that these differences are non-negative.

We will interchangeably think of a GT pattern as a semistandard Young tableau (SSYT). In this perspective, $(T^j_1, T^j_2, \cdots, T^j_j)$ is the partition formed by the cells of the tableau which contain entries $\leq j$. It follows that  $\NE_{ij}(T)$ is the number of cells in the $i^{\,\text{th}}$ row of the tableau which contain the entry $j+1$. We let $x^T$ denote the $x$-weight of the corresponding tableau. The following fermionic formula for the \qw polynomial appears in \cite{HKKOTY, Kirillov_newformula} and follows readily from Macdonald's more general formula \cite[Chap VI, (7.13)']{MacMainBook}:
\begin{equation}\label{eq:ferm}
  \wl = \sum_{T \in \GT} x^T \prod_{1 \leq i \leq j <n}\qbinom{NE_{ij}(T) + SE_{ij}(T)}{NE_{ij}(T)}
\end{equation}
Following \cite{karpthomas}, we define $\wt_q(T) = \displaystyle\prod_{1 \leq i \leq j <n}\displaystyle\qbinom{NE_{ij}(T) + SE_{ij}(T)}{NE_{ij}(T)}$.

\subsection{Partition overlaid patterns}\label{sec:pop-qbinom}
We recall that the $q$-binomial $\qbinom{k+\ell}{k}$ is the generating function of partitions that fit into a $k \times \ell$ rectangle, i.e.,   $\qbinom{k+\ell}{k} = \sum q^{|\gamma|}$ where $\gamma = (\gamma_1 \geq \gamma_2\geq \cdots \geq\gamma_k \geq 0)$ with $\ell \geq \gamma_1$. We also identify partitions of the above form with strictly decreasing $k$-tuples of integers between $0$ and $k+\ell-1$ via the bijection $\gamma \mapsto \bgam = \gamma+\delta$ where $\delta=(k-1, k-2, \cdots, 0)$. 

\begin{figure}
  \begin{center}
\begin{tikzpicture}[x={(1cm*0.5,-\rootthree cm*0.5)},y={(1cm*0.5,\rootthree cm*0.5)}]
  \foreach\i in{0,...,3}
\foreach\j in{\i,...,3}{
  \pgfmathtruncatemacro{\k}{4 - \j + \i};
  \pgfmathtruncatemacro{\l}{\i + 1};
  \draw(\i,\j)node(a\i\j){$\scriptstyle{T^{\k}_{\l}}$};
}
  \foreach\i/\ii in{0/-1,1/0,2/1,3/2}
 \foreach\j/\jj in{0/-1,1/0,2/1,3/2}{
  \ifnum\i<\j \draw[color=red,thick](a\i\jj)--(a\i\j); \fi
  \ifnum\i>\j\else\ifnum\i>0 \draw[color=blue](a\ii\j)--(a\i\j);\fi\fi
 }
\end{tikzpicture}\qquad
\begin{tikzpicture}[x={(1cm*0.5,-\rootthree cm*0.5)},y={(1cm*0.5,\rootthree cm*0.5)}]
  \draw(0,0)node(a00){$\scriptstyle{10}$};\draw(0,1)node(a01){$\scriptstyle{8}$};\draw(0,2)node(a02){$\scriptstyle{7}$};\draw(0,3)node(a03){$\scriptstyle{4}$};
  \draw(1,1)node(a11){$\scriptstyle{6}$};\draw(1,2)node(a12){$\scriptstyle{5}$};\draw(1,3)node(a13){$\scriptstyle{2}$};
  \draw(2,2)node(a22){$\scriptstyle{4}$};\draw(2,3)node(a23){$\scriptstyle{2}$};
  \draw(3,3)node(a33){$\scriptstyle{0}$};
  \foreach\i/\ii in{0/-1,1/0,2/1,3/2}
 \foreach\j/\jj in{0/-1,1/0,2/1,3/2}{
  \ifnum\i<\j \draw[color=red,thick](a\i\jj)--(a\i\j); \fi
  \ifnum\i>\j\else\ifnum\i>0 \draw[color=blue](a\ii\j)--(a\i\j);\fi\fi
}
\end{tikzpicture}
\qquad
\begin{tikzpicture}[x={(1cm*0.5,-\rootthree cm*0.5)},y={(1cm*0.5,\rootthree cm*0.5)}]
    \ytableausetup{mathmode, boxsize=0.5em}
  \draw(0,0)node {\ydiagram{1,1}};\draw(0,1)node{\ydiagram{2}};\draw(0,2)node{\ydiagram{2,1}};
  \draw(1,1)node {\ydiagram{1}};\draw(1,2)node{$\emptyset$};
  \draw(2,2)node {\ydiagram{2,2}};
\end{tikzpicture}
\caption{A GT pattern for $n=4$. The NE and SE differences are those along the red and blue lines. On the right is a partition overlay compatible with this GT pattern.}
\label{fig:gt-ssyt}
\end{center}
\end{figure}
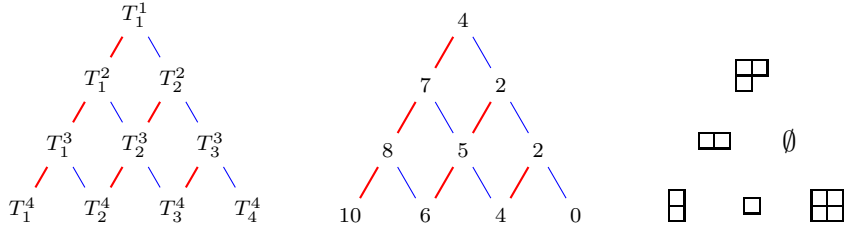

As shown in \cite{RRV-CLpop}, the right-hand side of \eqref{eq:ferm} can be interpreted in terms of the so-called {\em partition overlaid patterns} (POPs). A POP of shape $\lambda$ is a pair $(T, \Lambda)$ where $T \in \GT$ and $\Lambda=(\Lambda_{ij}: 1 \leq i \leq j <n)$ is a tuple of partitions such that each $\Lambda_{ij}$  fits into a rectangle of size $NE_{ij}(T) \times SE_{ij}(T)$. For example, if $T$ is the GT pattern of Figure~\ref{fig:gt-ssyt}, we could take $\Lambda_{11} = (2,1,0),\;  \Lambda_{12} = (2),\;  \Lambda_{13} =(1,1),\;  \Lambda_{22} = (0,0,0) ,\; \Lambda_{23} = (1), \; \Lambda_{33} = (2,2)$. We imagine the $\Lambda_{ij}$ as being placed in a triangular array as in Figure~\ref{fig:gt-ssyt}.
We let $\pop$ denote the set of POPs of shape $\lambda$.
It is now clear from \eqref{eq:ferm} that 
\begin{equation}\label{eq:fermpop}
  \wl = \sum_{(T, \Lambda) \in \pop} x^T q^{|\Lambda|}
  \end{equation}
where $|\Lambda| = \sum_{i,j} |\Lambda_{ij}|$.
We remark that $\wl$ is the character of the {\em local Weyl module} $W_{\mathrm{loc}}(\lambda)$ - a module for the current algebra $\mathfrak{sl}_n[t]$ \cite{ChariLoktev-original,ChariIon-BGG}. Further, POPs of shape $\lambda$ index a special basis of this module with Gelfand-Tsetlin like properties \cite{ChariLoktev-original,RRV-CLpop}.  

\subsection{Projection and Branching for Partition overlaid patterns}
Given $\lambda = (\lambda_1 \geq \lambda_2 \geq \cdots \geq \lambda_n \geq 0)$, we say that $\mu = (\mu_1, \mu_2, \cdots, \mu_{n-1})$ {\em interlaces} $\lambda$ (and write $\mu \prec \lambda$) if $\lambda_i \geq \mu_i \geq \lambda_{i+1}$ for $1 \leq i <n$. 
The \qw polynomials have the following important properties which readily follow from \eqref{eq:fermpop}:

\smallskip
\noindent
({\em projection}) $\wl[q=0] = s_\lambda(X_n)$, the Schur polynomial, and
\begin{equation}\label{eq:branch}
\!\!\text{({\em branching}) } W_\lambda(x_1, x_2, \cdots, x_{n-1}, x_n=1; q) = \sum_{\mu \prec \lambda} \prod_{1 \leq i <n} \qbinom{\lambda_i-\lambda_{i+1}}{\lambda_i-\mu_i} \cdot W_\mu(X_{n-1}; q)
\end{equation}
In fact, Chari-Loktev \cite[]{ChariLoktev-original} lift \eqref{eq:branch} to the level of modules, showing that the local Weyl module $W_{\mathrm{loc}}(\lambda)$ when restricted to $\mathfrak{sl}_{n-1}[t]$ admits a filtration whose successive quotients are of the form $W_{\mathrm{loc}}(\mu)$ for $\mu \prec \lambda$; further their graded multiplicities are precisely given by the product of $q$-binomial coefficients that appear in \eqref{eq:branch}.

The combinatorial shadow of projection is the map $\pr: \pop \to \GT$ given by $\pr(T,\Lambda) = T$. Likewise, we define {\em combinatorial branching} to be the map $\br: \pop \to \bigsqcup_{\mu \prec \lambda} \pop[\mu]$ defined by $\br(T,\Lambda) = (T^\dag, \Lambda^\dag)$ where $T^\dag$ is obtained from $T$ by deleting its bottom row, and $\Lambda^\dag$ is obtained from $\Lambda$ by deleting the overlays $\Lambda_{ij}$ with $j=n-1$.
\subsection{Box complementation} In addition to $\pr$ and $\br$, $\pop$ is endowed with another important map, which we term {\em box complementation}. Observe that given a partition $\pi = (\pi_1 \geq \pi_2 \geq \cdots \geq \pi_k \geq 0)$ fitting into a $k \times \ell$ rectangle, i.e., with $\pi_1 \leq \ell$, we may consider its complement in this rectangle, defined by $\pi^c = (\ell-\pi_k \geq \ell-\pi_{k-1} \geq \cdots \geq \ell-\pi_1)$. Now, for $(T,\Lambda) \in \pop$, define $\bcomp(T,\Lambda)=(T,\Lambda^c)$ where for each $i,j$, $(\Lambda^c)_{ij}$ is defined to be the complement of $\Lambda_{ij}$ in its bounding rectangle of size $\NE_{ij}(T) \times \SE_{ij}(T)$.

We note that since $|\Lambda|\neq |\Lambda^c|$ in general, $\bcomp$ preserves $x$-weights, but not $q$-weights. However $|\Lambda| + |\Lambda^c| = \sum_{i,j} \NE_{ij}(T)\SE_{ij}(T) =:\area(T)$ (in the terminology of \cite{RRV-CLpop}), which depends only on $T$.

\section{Projection and branching for Column strict fillings}
Our goal is to construct natural bijections between $\csf$ and $\pop$ which explain the equality of \eqref{eq:wlinvquinv} and \eqref{eq:fermpop} for $v=\inv, \quinv$. In addition to preserving $x$- and $q$-weights, we would like our bijections to be compatible with projection and branching. Towards this end, we first define these latter maps in the setting of $\csf$.

\subsection{Projection: rowsort} Given $F \in \csf$, let $\rsort(F)$ denote the filling obtained from $F$ by sorting entries of each row in ascending order. In light of the following easy lemma, we think of $\rsort$ as the projection map in the CSF setting.
\begin{lemma}
  If $F \in \csf$, then $\rsort(F) \in \SSYT(\lambda) \cong \GT$.
\end{lemma}

\subsection{Branching: delete-and-splice}\label{sec:dsplice}
A strictly increasing sequence $a=(a_1 < a_2 < \cdots < a_m)$ of positive integers will also be termed a {\em column tuple} with $\len(a)=m \geq 0$. Let $\ell  \geq 1$ and suppose $\sigma=(\sigma_1 < \sigma_2 < \cdots < \sigma_{\ell-1})$ and $\tau=(\tau_1 < \tau_2 < \cdots < \tau_{\ell})$ are column tuples of length $\ell-1$ and $\ell$ respectively. We set $\sigma_0=0$ and let $k$ denote the maximum element of the (non-empty) set $\{1 \leq i \leq \ell: \sigma_{i-1} < \tau_i\}$. Define $\splice(\sigma,\tau) = (\barsig,\bartau)$ where
\[ \barsig_i = \begin{cases} \sigma_i & 1\leq i <k \\ \tau_i &k \leq i \leq \ell \end{cases} \;\;\;\;\;\;\; \text{ and } \;\;\;\;\;\; \bartau_i = \begin{cases} \tau_i & 1\leq i <k \\ \sigma_i &k \leq i < \ell \end{cases}\]
i.e., $\barsig,\bartau$ are obtained by swapping certain suffix portions of $\sigma, \tau$. The choice of $k$ ensures that $\barsig, \bartau$ are also column tuples; we also have $\len(\barsig) = \len(\tau)$ and $\len(\bartau)=\len(\sigma)$.   \ytableausetup{mathmode, smalltableaux} For instance, when $(\sigma,\tau) = (\; \ytableaushort{1,5} \; ,\, \ytableaushort{2,3,4}\;)$, we get $(\barsig,\bartau) = (\; \ytableaushort{1,3,4}\; ,\, \ytableaushort{2,5}\;)$.

We now define the delete-and-splice rectification (``dsplice") map on $F \in \csf$ as follows:
{\bf (1)} delete all cells in $F$ containing the entry $n$ and let $\fdag$ denote the resulting filling. While its column entries remain strictly increasing, $\fdag$ may no longer be of partition shape. {\bf (2)} Let $\sigma^{(j)}\, (j \geq 1)$ denote the column tuple obtained by reading the $j^{\text{\,th}}$ column of $\fdag$ from top to bottom.  If $\fdag$ is not of partition shape, there exists $j \geq 1$ such that $\len(\sigma^{(j+1)}) = \len(\sigma^{(j)}) +1$. Choose any such $j$ and modify $\fdag$ by replacing the pair of columns $(\sigma^{(j)}, \sigma^{(j+1)})$ in $\fdag$ by $\splice(\sigma^{(j)}, \sigma^{(j+1)})$.  This swaps the column lengths and brings the shape of $\fdag$ one step closer to being a partition. {\bf (3)} If the shape of $\fdag$ is a partition, STOP. Else go back to step 2.

It is clear that this process terminates and finally produces a CSF of partition shape (filled by numbers between $1$ and $n-1$), which we denote $\dsplice(F)$. The following properties hold:

\begin{proposition}\label{prop:dsplice-prop}
With notation as above: (i) $D:=\dsplice(F)$ is independent of the intermediate choices of $j$ made in step 2 of the procedure. (ii) $\rsort(D)$ is obtained from $\rsort(F)$ by deleting the cells containing the entry $n$. (iii) If $\mu$ and $\lambda$ are the shapes of $D$ and $F$ respectively, then $\mu \prec \lambda$.
\end{proposition}

We consider $\dsplice$ to be the combinatorial branching map in the CSF context. Its key property is its compatibility with the natural branching map $\br$ of the POP setting (Theorem~\ref{thm:mainthm} below). While each $\splice$ operation is ``local'', the end result $\dsplice(F)$ can have a fair bit of ``intermixing'' amongst columns of $F$ (see also \S\ref{sec:concl-rem} for a pictorial description).

\section{The main theorem}

\begin{thm}\label{thm:mainthm}
  For any $n \geq 1$ and any partition $\lambda: \lambda_1 \geq \lambda_2 \geq \cdots \geq \lambda_n \geq 0$ with at most $n$ nonzero parts, there exist two bijections $\psi_{\inv}$ and $\psi_{\quinv}$ from $\csf$ to $\pop$ with the following properties:

  \medskip
  \noindent
    1. If $\psi_v(F) = (T,\Lambda)$, then $x^F = x^T$ and $v(F) = |\Lambda|$, for $v=\inv$ or $\quinv$.

\medskip
\noindent
2. The following diagrams commute ($v=\inv$ or $\quinv$):
\begin{enumerate}
  \item[(A)]
\[\begin{tikzcd}
	\csf && \pop \\
	& \GT
	\arrow["{\psi_v}", from=1-1, to=1-3]
	\arrow["\rsort"'{pos=0.4}, from=1-1, to=2-2]
	\arrow["\pr"{pos=0.4}, from=1-3, to=2-2]
\end{tikzcd}\]

\item[(B)] \[\begin{tikzcd}
	\csf &&& \pop \\
	\\
	\displaystyle\bigsqcup_{\mu \prec \lambda}\csf[\mu] &&& \displaystyle\bigsqcup_{\mu \prec \lambda}\pop[\mu]
	\arrow["{\psi_v}", from=1-1, to=1-4]
	\arrow["\dsplice"', from=1-1, to=3-1]
	\arrow["{\psi_v}"', from=3-1, to=3-4]
	\arrow["\br", from=1-4, to=3-4]
\end{tikzcd}\]
\end{enumerate}

\medskip
\noindent
3. The two bijections are related via the commutative diagram:
\[\begin{tikzcd}
	& \csf \\
	\pop && \pop
	\arrow["{\psi_{\mathrm{inv}}}"', from=1-2, to=2-1]
	\arrow["{\psi_{\mathrm{quinv}}}", from=1-2, to=2-3]
	\arrow["\bcomp", from=2-1, to=2-3]
\end{tikzcd}\]\qed
\end{thm}

To summarize, $\psi_{\mathrm{inv}}$ and $\psi_{\mathrm{quinv}}$ acting on a CSF produce POPs with the same underlying GT pattern, but with complementary overlays. These bijections are compatible with the natural projection and branching maps, and preserve $x$- and appropriate $q$-weights (inv or quinv). Note the slight abuse of notation in part 2(B) above: for $\mu \prec \lambda$, $\csf[\mu]$ denotes the set of column strict fillings $F: \dg(\mu)\to [n-1]$ (rather than $[n]$). Theorem~\ref{thm:mainthm}, with the exception of part 2(B), can also be formulated in the setting of $q$-Whittaker functions in infinitely many variables.
Next, we obtain the following corollaries:
\begin{corollary}
Let $T \in \GT$ and let $\rsort^{-1}(T) = \{F \in \csf: \rsort(F)=T\}$ be the fiber of $\rsort$ over $T$. 
\begin{enumerate}
\item $\displaystyle\sum_{F \in \rsort^{-1}(T)} q^{\inv(F)} = \displaystyle\sum_{F \in \rsort^{-1}(T)} q^{\quinv(F)} = \wt_q(T)$.
\item $\inv(F) + \quinv(F) = \area(T)$ is  constant for $F \in \rsort^{-1}(T)$.
\end{enumerate}
\end{corollary}
An interpretation of $\wt_q(T)$ in terms of flags of subspaces compatible with nilpotent operators appears in \cite[Theorem 5.8(i)]{karpthomas}.
In \cite{AMM}, the authors asked for an explicit bijection on $\scf$ which interchanges the $\inv$ and $\quinv$ statistics. We describe this bijection on $\csf$, thereby partially answering their question.
\begin{corollary}
The map $\Omega: \psi_{inv}^{-1} \circ \psi_{\quinv} = \psi_{inv}^{-1} \circ \, \bcomp \circ \, \psi_{\inv}: \csf \to \csf$ is an involution satisfying $\inv(\Omega(F)) = \quinv(F)$ for all $F \in \csf$.
\end{corollary}
The explicit construction of the $\psi_v$ and their inverses in the next section makes $\Omega$ effectively computable.

\section{Proof sketch}
\subsection{Cellwise zcounts and quinv triples}
\begin{figure}
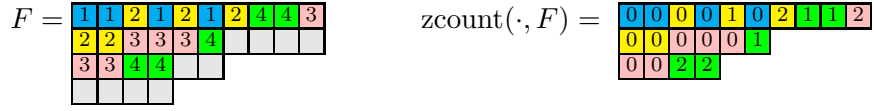

  \begin{center}
  \ytableausetup{mathmode, smalltableaux}
$F=\;$\begin{ytableau}
  *(cyan) 1  & *(cyan) 1 & *(yellow) 2 & *(cyan) 1 & *(yellow) 2 & *(cyan) 1 & *(yellow) 2 & *(green) 4 & *(green) 4 & *(pink) 3\\
  *(yellow) 2  & *(yellow) 2 & *(pink) 3 & *(pink) 3 & *(pink) 3 & *(green) 4 &*(light-gray) &*(light-gray) &*(light-gray) &*(light-gray) \\
  *(pink) 3  & *(pink) 3 & *(green) 4 & *(green) 4 &*(light-gray) &*(light-gray) \\
  *(light-gray) &*(light-gray) &*(light-gray) &*(light-gray) 
  \end{ytableau} \qquad \quad
 $\zcount(\cdot,F)=\;$ \begin{ytableau}
  *(cyan) 0  & *(cyan) 0 & *(yellow) 0 & *(cyan) 0 & *(yellow) 1 & *(cyan) 0 & *(yellow) 2 & *(green) 1 & *(green) 1 & *(pink) 2\\
  *(yellow) 0  & *(yellow) 0 & *(pink) 0 & *(pink) 0 & *(pink) 0 & *(green) 1\\
  *(pink) 0  & *(pink) 0 & *(green) 2 & *(green) 2
\end{ytableau}
\caption{Here $F \in \csf$ for $\lambda=(10,6,4,0)$ and $n=4$. Cells of $F$ are coloured according to their entries. The gray cells are the extra cells in the augmented diagram $\dgh(\lambda)$. On the right are cellwise zcount values. Here $\quinv(F)=12$.}
  \label{fig:zcfig}
\end{center}
\end{figure}

We first describe the construction of $\psi_{\quinv}$.
For a partition $\lambda$, the augmented diagram $\dgh(\lambda)$ is $\dg(\lambda)$ together with one additional cell below the last cell in each column (see Figure~\ref{fig:zcfig}).
Given $F \in \csf$, a {\em quinv-triple} in $F$ is a triple of cells $(x,y,z)$ in $\dgh(\lambda)$ such that (i) $x, z \in \dg(\lambda)$ and $z$ is to the right of $x$ in the same row, (ii) $y$ is the cell immediately below $x$ in its column, (iii) $F(x) < F(z) < F(y)$, where we set $F(y) = \infty$ if $y$ lies outside $\dg(\lambda)$. It is easy to see that the quinv-triples considered in \cite{AMM} for $F \in \scf$ reduce to this description when $F$ is a CSF rather than a general filling. Thus,
$\quinv(F)$ as defined in \cite{AMM} equals the number of quinv-triples in $F$ (as defined above) for a CSF $F$.

Given $F \in \csf$, we define a function {\em zcount} which tracks the contributions of individual cells of $\dg(\lambda)$ to $\quinv(F)$ as follows: for each cell $c \in \dg(\lambda)$, let $\zcount(c, F) =$ the number of quinv-triples $(x,y,z)$ in $F$ with $z=c$. Clearly
\begin{equation}\label{eq:qwtquinv}
  \sum_{c \in \dg(\lambda)} \zcount(c,F) = \quinv(F)
\end{equation}

\noindent
We next group cells of the filling $F$ row-wise according to the entries they contain. More precisely,  let
$\cells(i,j,F) = \{c \in \dg(\lambda): c \text{ is in the } i^{th} \text{ row and } F(c) = j+1\}$ for $1 \leq i \leq j+1 \leq n$.
Figure~\ref{fig:zcfig} shows an example, with these groups colour-coded in each row.
It readily follows from \S\ref{sec:fermfor} that
\begin{equation}\label{eq:cellcard}
  |\cells(i,j,F)| = NE_{ij}(T),  \text{ where } T = \rsort(F).
\end{equation}
The next proposition brings the SE differences also into play \cite{brv-fullversion}:
\begin{proposition}\label{prop:zcounts}
  Let $F \in \csf$ and $T = \rsort(F)$. Fix $1 \leq i \leq j+1 \leq n$.
  \begin{enumerate}
  \item If $c \in \cells(i,j,F)$, then $\zcount(c,F) \leq \SE_{ij}(T)$.
  \item If $c, d \in  \cells(i,j,F)$ with $c$ lying to the right of $d$, then $\zcount(c,F) \geq \zcount(d,F)$.
    \item Further, equality holds in (1) for all $i, j$ and all cells $c \in \cells(i,j,F)$ iff $F=T$.
  \end{enumerate}
\end{proposition}

  \subsection{Definition of $\psi_{\quinv}$} We now have all the ingredients in place to define $\psi_{\quinv}$. Let $F \in \csf$ and $T =\rsort(F)$. For each $1 \leq i \leq j+1 \leq n$, consider the sequence
  \begin{equation}\label{eq:lamij}
    \Lambda_{ij} = (\zcount(c,F): \; c \in \cells(i,j,F) \text{ traversed right to left in row } i).
    \end{equation}
In Figure~\ref{fig:zcfig}, this amounts to reading the entries of a fixed colour from right to left in a given row of $\zcount(\cdot, F)$.
  By Proposition~\ref{prop:zcounts}, this is a weakly decreasing sequence bounded above by $\SE_{ij}(T)$. Together with \eqref{eq:cellcard}, this implies that $\Lambda_{ij}$ may be viewed as a partition fitting into the $\NE_{ij}(T) \times \SE_{ij}(T)$ rectangle. Since $\SE_{ij} =0$ for $i=j+1$, $\Lambda_{ij}$ is the zero sequence in this case. We drop the pairs $(j+1,j)$ to conclude that if $\Lambda = (\Lambda_{ij}: 1 \leq i \leq j <n)$, then $(T, \Lambda) \in \pop$. We define $\psi_{quinv}(F) = (T, \Lambda)$. 
  Clearly, $x^F = x^T$ and \eqref{eq:qwtquinv} implies $\quinv(F) = |\Lambda|$, establishing (1) of Theorem~\ref{thm:mainthm} for $v=\quinv$.
  
\subsection{refinv triples}
\begin{figure}
  \begin{center}
    \begin{tikzpicture}
      \ytableausetup{mathmode, smalltableaux}
      \draw(0,0)node {\ydiagram[*(yellow)]{1,1}};
      \draw(2,0.2)node{\ydiagram[*(yellow)]{1}};
      \draw[dashed](0.3,0.2)--(1.7,0.2);
    \end{tikzpicture}
    \qquad
        \begin{tikzpicture}
      \ytableausetup{mathmode, smalltableaux}
      \draw(0,0)node {\ydiagram[*(cyan)]{1,1}};
      \draw(2,-0.2)node{\ydiagram[*(cyan)]{1}};
      \draw[dashed](0.3,-0.2)--(1.7,-0.2);
    \end{tikzpicture}
\qquad
        \begin{tikzpicture}
      \ytableausetup{mathmode, smalltableaux}
      \draw(2,0)node {\ydiagram[*(orange)]{1,1}};
      \draw(0,0.2)node{\ydiagram[*(orange)]{1}};
      \draw[dashed](0.3,0.2)--(1.7,0.2);
    \end{tikzpicture}
        \caption{(left to right) Configuration of quinv, inv and refinv triples. }
        \label{fig:triplesconfig}
  \end{center}
\end{figure}
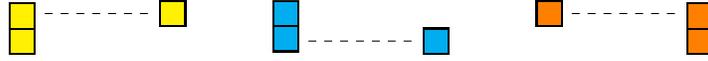

We now turn to the definition of $\psi_{\inv}$. While we may anticipate doing this via a modification of the foregoing arguments, replacing quinv-triples with Haglund-Haiman-Loehr's inv-triples, that turns out not to work out-of-the-box. In place of the latter (see Figure~\ref{fig:triplesconfig}), we consider triples $(x,y,z)$ in $\dgh(\lambda)$ where (i) $x, z \in \dg(\lambda)$ with $z$ to the left of $x$ in the same row, (ii) $y$ is the cell immediately below $x$ in its column. Given $F \in \csf$, we call $(x,y,z)$ a {\em refinv-triple} (or ``reflected inv-triple'') for $F$ if in addition to (i) and (ii), we also have (iii) $F(x) < F(z) < F(y)$, where $F(y):=\infty$ if $y \not\in \dg(\lambda)$. We have \cite{brv-fullversion}:
\begin{proposition}\label{prop:refinv}
  For $F \in \csf$, $\inv(F)$ equals the number of refinv-triples of $F$.
\end{proposition}

\begin{remarks}
  1. We may in fact define a new statistic\footnote{In fact, {\em refquinv} can also be likewise defined on all fillings, and agrees with {\em quinv} on CSFs. But rephrased in terms of refquinv-triples, this involves counting such triples with signs \cite{brv-fullversion}.}  {\em refinv} on all fillings $F \in \scf$ as follows: $\refinv(F) = \Inv(F) - \sum_{u \in \Des F} \coarm(u)$, borrowing notation of \cite[\S 2]{HHL-I}. This replaces {\em arm} in HHL's definition by {\em coarm}. The content of Proposition~\ref{prop:refinv} is that $\refinv(F)=\inv(F)$ for $F \in \csf$. In fact, this equality holds more generally for all fillings $F$ whose descent set is a union of rows of $\dg(\lambda)$. 
  More generally, The proof follows from the observation that the sum of $\arm(c)$ equals the sum of $\coarm(c)$ when the cell $c$ runs over the union of a subset of rows of $\dg(\lambda)$.

  \medskip
  \noindent
  2. The refinv triples for $F \in \csf$ actually make an appearance in \cite[\S 2.2]{Kirillov_newformula}, where they are attributed to Zelevinsky (and their total number denoted $\widetilde{ZEL}$). From this perspective, the content of Proposition~\ref{prop:refinv} is that $\widetilde{ZEL}(F) = \inv(F)$.
\end{remarks}

\subsection{$\zcb$, $\zcount$ and the proof of the main theorem}
Given $F \in \csf$ and $c \in \dg(\lambda)$, define  $\zcb(c, F) =$ the number of refinv-triples $(x,y,z)$ in $F$ with $z=c$. In light of Proposition~\ref{prop:refinv}, it is clear that
\begin{equation}\label{eq:qwtinv}
\sum_{c \in \dg(\lambda)} \zcb(c,F) = \inv(F)
\end{equation}
We have the following relation between $\zcb$ and $\zcount$ \cite{brv-fullversion}:
\begin{proposition}\label{prop:zcandzcb}
Let $F \in \csf$ and $T=\rsort(F)$. Let $1 \leq i \leq j+1 \leq n$ and $c \in \cells(i,j,F)$. Then $\zcount(c,F) + \zcb(c,F) = \SE_{ij}(T)$.
\end{proposition}

We may now define $\psi_{\inv}$ following the template of $\psi_{\quinv}$. Given $F \in \csf$, let $T =\rsort(F)$. For each $1 \leq i \leq j < n$,  consider the sequence:
  \[\overline{\Lambda}_{ij} = (\zcb(c,F): \; c \in \cells(i,j,F) \text{ traversed left to right in row } i)\]
   Recall also the definition of the partition $\Lambda_{ij}$ from \eqref{eq:lamij}. It follows from Propositions~\ref{prop:zcounts} and \ref{prop:zcandzcb} that $\overline{\Lambda}_{ij}$ is the box-complement of $\Lambda_{ij}$ in the $\NE_{ij}(T) \times \SE_{ij}(T)$ rectangle. Letting $\overline{\Lambda} = (\overline{\Lambda}_{ij}: 1 \leq i \leq j<n)$, we define $\psi_{\inv}(F) = (T,\overline{\Lambda})$. As in the case of $\quinv$, we have $x^F = x^T$, and  $\inv(F) = |\overline{\Lambda}|$ by \eqref{eq:qwtinv}. This proves part (1) of Theorem~\ref{thm:mainthm} for $v=\inv$.

Since by definition $\pr (\psi_v(F)) = T$ for $v=\inv, \quinv$,  Part (2A) of Theorem~\ref{thm:mainthm} follows. Part (3) of Theorem~\ref{thm:mainthm}  follows from the fact that $\Lambda$ and $\overline{\Lambda}$ are box complements of each other in the appropriate rectangles.
That the diagrams in part (2B) of Theorem~\ref{thm:mainthm} are commutative follows from an analysis of each elementary splice step of the $\dsplice$ map; we defer the details to \cite{brv-fullversion}.

Finally, this leaves us with proving that the $\psi_v$ are bijections. We sketch the construction of $\psi_{\inv}^{-1}$. Given $(T, \Lambda) \in \pop$, construct the filling $F:=\psi_{\inv}^{-1}(T,\Lambda) \in \csf$ inductively row-by-row, from the bottom ($n^{th}$) row to the top as follows: (a) fill all cells of the $n^{th}$ row (if nonempty) with $n$, (b) let $1 \leq i \leq j < n$; assuming that all rows of $F$ strictly below row $i$ have been completely determined and that the locations of entries $>(j+1)$ in row $i$ have been determined, we now need to fill $\NE_{ij}(T)$ many cells of row $i$ with the entry $j+1$. It turns out that the number of cells in row $i$ in which we can potentially put a $j+1$ without violating the CSF condition thus far is exactly $k+\ell$ where $k=\NE_{ij}(T)$ and $\ell=\SE_{ij}(T)$. We label these cells $0, 1, \cdots, k+\ell-1$ from right to left (left-to-right when defining $\psi_{\quinv}^{-1}$). We now use the identification from \S\ref{sec:pop-qbinom} of partitions fitting inside a $(k \times \ell)$-box with $k$-tuples of distinct integers in $0, 1, \cdots, k+\ell-1$. Via this, the partition $\Lambda_{ij}$ can be viewed as a $k$-tuple of candidate cells in row $i$; we put the entry $j+1$ into these, (c) fill the remaining cells of row $i$ with the entry $i$. The rest of the argument is straightforward \cite{brv-fullversion}. \qed

For example, let $n=4$, $\lambda = (10, 6, 4,0)$ and let $T, \Lambda$ be the GT pattern and overlay depicted in Figure~\ref{fig:gt-ssyt}. 
Then $\psi_{\quinv}^{-1}(\mathcal{T},\Lambda)$ is precisely the CSF $F$ of Figure~\ref{fig:zcfig}, while

\smallskip
\noindent
\ytableausetup{mathmode, smalltableaux}
$\psi_{\inv}^{-1}(\mathcal{T},\Lambda) = \ytableaushort{2111321442,332243,4433}$ 

\section{Local Weyl modules and limit constructions}
Finally, we can apply these ideas to the study of local Weyl modules, in particular to the {\em limit constructions} of \cite{fourier-littelmann,RRV-CLpop,ravinder2018stability}. Let $L(\Lambda_0)$ denote the basic representation of the affine Lie algebra $\widehat{\mathfrak{sl}_n}$ \cite[Prop. 12.13]{kac}. Using Theorem~\ref{thm:mainthm} to replace POPs with CSFs as our model in \cite[Corollary 5.13]{RRV-CLpop}, we deduce \cite{brv-fullversion}:
\begin{proposition}
  Fix $n \geq 2$ and consider the partition $\theta=(2,1,1,\cdots,1,0)$ with $n-1$ nonzero parts and $|\theta|=n$. For $k \geq 0$, let $\mathcal{C}_k$ denote the set of CSFs $\,F$ of shape $k\theta$ and entries in $[n]$, with the property that either $1$ occurs in the first column of $F$ or $1$ does not occur in its last column. Then
    $\sum_{k \geq 0} \sum_{F \in \mathcal{C}_k} \,x^F \,q^{k^2-\inv(F)}$ equals the character of $L(\Lambda_0)$.
\end{proposition}
There is also a more general version with $\lambda+k\theta$ in place of $k\theta$ (for appropriate $\lambda$), mirroring \cite[Corollary 5.13]{RRV-CLpop}.

\section{Concluding Remarks}\label{sec:concl-rem}
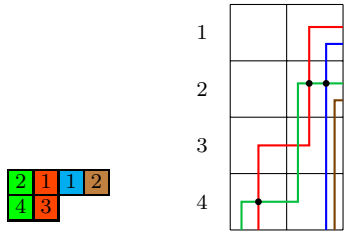
\begin{figure}
  \begin{center}
    \begin{tikzpicture}
      \draw(0,0) node{\begin{ytableau} *(green)2 & *(orangered)1 & *(cyan)1 & *(brown)2\\*(green)4 & *(orangered)3\end{ytableau}};
    \end{tikzpicture}
    \qquad
    \begin{tikzpicture}[scale=0.75]
      \draw[] (0,0) grid (2,4);
      \draw[color=red,thick](2,3.6)--(1.4,3.6)--(1.4,1.5)--(0.5,1.5)--((0.5,0);
      \draw[color=blue,thick](2,3.3)--(1.7,3.3)--(1.7,0);
      \draw[color=darkpastelgreen,thick](2,2.6)--(1.2,2.6)--(1.2,0.5)--(0.2,0.5)--((0.2,0);
      \draw[color=chocolate,thick](2,2.3)--(1.85,2.3)--(1.85,0);
      \fill (1.4,2.6) circle [radius=0.6mm];
      \fill (1.7,2.6) circle [radius=0.6mm];
      \fill (0.5,0.5) circle [radius=0.6mm];
      \foreach\i in{0,...,3}{
        \pgfmathtruncatemacro{\k}{4 - \i};
        \draw(-0.5,\i+0.5)node {$\scriptstyle \k$};};
    \end{tikzpicture}
    \caption{A CSF $F$ with columns colour-coded to match its lattice path representation. The three marked intersections show that $\inv(F)=3$.}
    \label{fig:csf-wheelerdig}
  \end{center}
  \end{figure}

For the modified Hall-Littlewood polynomials $\qpl[\lambda']$ of \eqref{eq:qlpinvquinv}, the fermionic formula appears in \cite[(0.2)]{Kirillov_newformula}. Analogous to \eqref{eq:fermpop}, this can now be recast as a {\em weighted sum} over {\em partition overlaid plane-partitions} (POPP) of shape $\lambda$. Theorem~\ref{thm:mainthm} takes the form of bijections from $\wdf$ to $\popp$ (or equivalently, from tabloids to partition overlaid reverse-plane-partitions). The subtlety here is that POPPs need to be weighted with an additional power of $q$ (which depends only on the underlying plane-partition, cf \cite[(0.2)]{Kirillov_newformula}). The refinv- or quinv-triples in this case also involve $\leq$ relations (rather than just $<$) and this extra $q$-power keeps track of certain equalities among the triples \cite{brv-fullversion}.

Secondly, the bijections of Theorem~\ref{thm:mainthm} (and those indicated above for the modified Hall-Littlewood case) have an attractive interpretation in terms of lattice-path diagrams \cite{GarbaliWheeler,borodin-wheeler}.
Figure~\ref{fig:csf-wheelerdig} shows the lattice path representation of a CSF $F$; $\inv(F)$ is just the total number of intersections of the form $\;
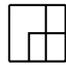
\begin{tikzpicture}[scale=0.75] \draw[thick] (0,0) grid (1,1); \draw[thick](0.65,0)--(0.65,1); \draw[thick](1,0.5)--(0.35,0.5)--(0.35,0); \end{tikzpicture}\;$ in the grid, and refining this further to each box of the grid produces the partition overlay as well \cite{brv-fullversion}. Likewise $\quinv(F)$ counts non-intersections of the above form. The $\dsplice$ map of \S\ref{sec:dsplice} translates into deletion of the last row of the grid followed by appropriate rectifications  $\begin{tikzpicture}[scale=0.75] \draw[thick] (0,0) grid (1,1); \draw[thick,color=red](0,0.5)--(0.75,0.5)--(0.75,1); \draw[thick,color=darkpastelgreen](0.5,1)--(0.5,0); \draw[thick,->](1.2,0.5)--(2.2,0.5); \end{tikzpicture}\; 	 \begin{tikzpicture}[scale=0.75] \draw[thick] (0,0) grid (1,1); \draw[thick,color=red](0,0.5)--(0.5,0.5)--(0.5,1); \draw[thick,color=darkpastelgreen](0.75,1)--(0.75,0); \end{tikzpicture}$

\vspace{-4mm}
\printbibliography

\end{document}